\documentclass[11pt]{amsart}

\usepackage{amsmath}
\usepackage{amsfonts}
\usepackage{amssymb}
\usepackage{amsthm}
\usepackage{indentfirst}
\usepackage[T1]{fontenc}

\usepackage{tikz}
\usepackage{verbatim}
\usetikzlibrary{matrix}

\usepackage{enumerate}

\newtheorem{thm}{Theorem}

\newtheorem{conj}{Conjecture}

\newcommand{\bba}{{\mathbb{A}}}
\newcommand{\bbf}{{\mathbb{F}}}
\newcommand{\bbp}{{\mathbb{P}}}
\newcommand{\bbq}{{\mathbb{Q}}}

\newcommand{\bbz}{{\mathbb{Z}}}
\newcommand{\bbc}{{\mathbb{C}}}
\newcommand{\Gal}{\operatorname{Gal}}
\newcommand{\Aut}{\operatorname{Aut}}
\newcommand{\GL}{\operatorname{GL}}

\newcommand{\rk}{\operatorname{rk}}

\newcommand{\bs}{\boldsymbol{\sigma}}

\newcommand{\Hom}{\operatorname{Hom}}
\newcommand{\Span}{\operatorname{Span}}
\newcommand{\tor}{\operatorname{tor}}
\newcommand{\ord}{\operatorname{ord}}

\title{Abelian Varieties and Finitely Generated Galois Groups}

\author{Bo-Hae Im}
\address{Department of Mathematical Sciences, KAIST, 291 Daehak-ro, Yuseong-gu, Daejeon, 34141, South Korea}
\email{bhim@kaist.ac.kr}

\author{Michael Larsen}
\address{Department of Mathematics, Indiana University, Bloomington, IN, 47405, U.S.A.}
\email{mjlarsen@indiana.edu}

\date{\today}

\begin{document}

\maketitle
\begin{abstract}
This paper surveys the methods that have been used to attack the conjecture, still open, that an abelian variety over a characteristic $0$ field with finitely generated Galois group is always of infinite rank.

\end{abstract}

\section{Introduction}

Let $K$ be a field finitely generated over $\bbq$.  By N\'eron's extension \cite[Chapter 6,~Theorem 1]{Lang} of the Mordell-Weil theorem, for every abelian variety $A/K$, the group $A(K)$ is finitely generated.

We say a field $K$ is \emph{anti-Mordell-Weil} (AMW) if for every non-trivial abelian variety $A/K$, the group
$A(K)$ is of infinite rank.

In 1974,  Frey and  Jarden \cite{FJ} proved the following theorem:

\begin{thm}
\label{prehistory}
If $K$ is finitely generated over $\bbq$ and $n$ is a positive integer, the set of $\bs := (\sigma_1,\ldots,\sigma_n)\in G_K^n$ such that $\bar K(\bs)$ is AMW is of measure $1$.
\end{thm}
Here $G_K := \Gal(\bar K/K)$, and $\bar K(\bs)$ denotes the fixed field of $\bar K$ under $\langle\bs\rangle := \langle \sigma_1,\ldots,\sigma_n\rangle$.

In \cite{Larsen}, one of us conjectured that more is true: that in fact $\bar K(\bs)$ is AMW for all $\bs$.  This conjecture remains open.
In this paper, we will present several different approaches to the problem, review what is known, and
discuss some variants and related open problems.

This conjecture can be reformulated as follows:

\begin{conj}
\label{conji}
Every characteristic $0$ field with finite transcendence degree and topologically finitely generated absolute Galois group is AMW.
\end{conj}

Indeed, $\bar K(\bs)$ has a finitely generated Galois group, and every characteristic $0$ field with finite transcendence degree and finitely generated absolute Galois group can be written in this form.

This implies the following, apparently stronger, statement:

\begin{conj}
\label{conj-ii}
Every characteristic $0$ field with topologically finitely generated absolute Galois group is AMW.
\end{conj}

Indeed, let $K$ be any characteristic $0$ field and $K' := \bar K(\bs)$ an algebraic extension of $K$ with a finitely generated Galois group.
Any abelian variety over $K'$ is obtained by extension of scalars from an abelian variety over a finitely generated subfield $k\subset K'$.
Let $k'$ and $\bar k$ denote the algebraic closure of $k$ in $K'$ and $\bar K$ respectively.
Then $\bar k$ is an algebraic closure of $k'$, and there is a natural continuous homomorphism $\pi\colon G_{K'}\to G_{k'}$.   As $\bar k$ is linearly
disjoint with $K'$ over $k'$ and algebraic over $k'$, we have $K'\otimes_{k'} \bar k\cong K'\bar k$ is a subfield of $\bar K$, and every $k'$-automorphism of $\bar k$ extends
uniquely to an automorphism of $K'\bar k$ which is trivial on $K'$ and thus to an automorphism of $\bar K$ which is trivial on $K'$.  We conclude that
$\pi$ is surjective, so $G_{k'}$ is finitely generated.  Thus Conjecture~\ref{conji} for $k'$ implies Conjecture~\ref{conj-ii}.

The following special case of Conjecture~\ref{conji} remains open, even in the case that $K=\bbq$ and $A$ is an elliptic curve:

\begin{conj}
\label{weak}
Let $A$ be an abelian variety over a number field $K$.  Then the rank of $A$ is infinite over $\bar K(\bs)$ for all $\bs$.
\end{conj}

There is also a characteristic $p$ version, but some care is needed because a locally finite field (i.e., a subfield of $\bar\bbf_p$) cannot be AMW.  For characteristic $p$ fields, one can define a field $K$ to be AMW if every abelian variety $A/K$ \emph{which is not isotrivial} has infinite rank over $K$.  There are no such abelian varieties over locally finite fields, so the condition is vacuous, and we can formulate:

\begin{conj}
\label{max}
Every field with topologically finitely generated absolute Galois group is AMW.
\end{conj}
Here the absolute Galois group of a field $K$ means $\Gal(\bar K/K)$, where $\bar K$ is a \emph{separable} closure of $K$.

The purpose of this paper is to survey what is known about these conjectures and related questions.  The emphasis will be on the wide variety of different methods which have been brought to bear, with varying degrees of success, on special cases.

\section{Probabilistic Methods}

We recall that a field $K$ is \emph{PAC} if every geometrically irreducible variety $V$ over $K$ has a $K$-point.  This implies that $V(K)$ is dense in $V$. We say $K$  is \emph{ample} (the terms \emph{large} and \emph{anti-mordellic} also appear in the literature) if every non-singular curve has either no $K$-points or infinitely many.  Any $K$-point on a curve which is irreducible over $K$ but not over some algebraic extension of $K$ must lie in at least two geometric components and therefore must be a singular point.  In other words, a $K$-point on any non-singular curve must lie on a geometrically connected irreducible component of that curve.  Thus PAC implies ample.

The ``PAC Nullstellensatz'' \cite[Theorem~18.6]{FriedJ} asserts that if $K$ is a countable field satisfying Hilbert approximation (e.g., any global field \cite[Chapter 9,~Theorem~4.2]{Lang}),
then with probability $1$, $\bar K(\bs)$ is PAC.
It follows that $\bar K(\bs)$ is ample.
Markus Junker and Jochen Koenigsmann \cite{JK}, made the following conjecture:

\begin{conj}
\label{ample}
Every infinite field with finitely generated Galois group is ample.
\end{conj}

Arno Fehm and Sebastian Petersen proved \cite{FP} that every ample field is AMW.  This shows that Conjecture~\ref{ample} implies Conjecture~\ref{max} and also gives a new proof of
Theorem~\ref{prehistory}.

The same approach can be used to give a probabilistic analysis of Galois representations coming from the non-torsion part of a non-trivial abelian variety $A$ over a field $K$ finitely generated over $\bbq$.
Let $V_A := A(\bar K)\otimes\bbq$, regarded as a space with discrete topology, and consider the
representation $\rho_A\colon G_K\to\Aut_{\bbq}V_A$.  This representation is continuous since
$$V_A := \sum_L A(L)\otimes \bbq = \bigcup_L A(L)\otimes \bbq,$$
where the union (or sum) is taken over all of the (countably many) finite Galois extensions $L/K$.
By N\'eron's theorem, each summand $A(L)\otimes \bbq$ is finite-dimensional.

Every ordered $n$-tuple $\bs\in G_K^n$ defines a homomorphism $e_{\bs}\colon F_n\to G_K$,
where $F_n$ is the free group on $n$ generators.  The composition $\rho_A\circ e_{\bs}$ is then a countable sum of finite dimensional $\bbq$-representations of $F_n$, each of which factors through a finite quotient
of $F_n$.  It is therefore a direct sum of irreducible $\bbq$-representations of $F_n$ which factor through finite quotients. We call such an irreducible representation an \emph{atom}.

By the \emph{generic representation} of $F_n$
we mean the direct sum of a countably infinite number of copies of every atom.  Thus,
$\rho_A\circ e_{\bs}$ is always a subrepresentation of the generic representation of $F_n$.

\begin{thm}
\label{new}
If $K$ is finitely generated over $\bbq$ and $n$ is a positive integer, the set of $\bs\in G_K^n$ for which the $G_K$-representation $\rho_A\circ e_{\bs}$ is generic for all non-trivial $A/K$ has measure $1$.
\end{thm}

Note that the space of $G_{\bar K(\bs)}$-invariants of $V_A$ is exactly $A(\bar K(\bs))\otimes\bbq$, so the infinite multiplicity of the trivial atom in $\rho_A\circ e_{\bs}$ is equivalent to the infinite rank of $A$ over $\bar K(\bs)$.

\begin{proof}
Since there are countably many isomorphism classes of abelian varieties over $K$, it suffices to prove the statement for a single one.  Likewise, it suffices to prove that each atom has infinite multiplicity in $V_A$ with probability $1$.  We therefore fix an atom $\alpha\colon F_n\to \GL(W)$ which factors through some finite quotient $G$ of $F_n$.  In particular, $G$ is generated by some $n$-element subset.

We fix an embedding $\iota\colon G\hookrightarrow A_N$, $N\ge 4$,
such that the restriction of the (unique) irreducible $N-1$-dimensional $\bbq$-representation of $A_N$ to $G$ contains $W$ as a subrepresentation.  This is possible for all $N\ge |G|+2$, by composing the regular permutation representation of $G$ with an embedding $S_{|G|}\hookrightarrow A_N$.
Note that $S_m$ embeds in $A_N$ for all $N\ge m+2\ge 3$,
and the permutation representation of $S_m$ is a subrepresentation of the restriction of the irreducible $N-1$-dimensional representation of $A_N$.

Since the $K$-points in any projective space are  Zariski-dense, Bertini's theorem implies that there exists a non-singular curve $X$ on $A$ defined over $K$ and passing through the identity $0$.
Let $g$ be the genus of $X$.
Since there are no rational curves on an abelian variety, we have $g\ge 1$.
We fix a non-empty $G_K$-stable finite subset $S$ of points in $X(\bar K)\setminus \{0\}$.

For any  $N \ge 2g+2|S|$, we consider the space of meromorphic functions $f$ on $X$ which vanish to order $\ge 2$ on every point of $S$, which are holomorphic except at $0$, and which has a pole of order $\le N$ at $0$.
By the Riemann-Roch theorem, the projective space of such functions has dimension $N-2|S|-g$; imposing the condition that $\ord_0 f\ge 1-N$ or the condition that $\ord_s f\ge  3$ for any specified $s\in S$ defines a projective subspace of codimension $1$, while imposing the condition that $\ord_t f\ge 2$ for some  $t\in X\setminus  (S\cup \{0\})$ defines a projective subspace of codimension $2$.
The union over $t$ of these codimension $2$ subspaces is therefore a subvariety of codimension $1$.

It follows that there is a dense open subvariety of functions $f$ on $X$
which have  a pole of exact order $N$ at $0$, no other poles, zeroes of exact order $2$ at each point in $S$, and no other multiple zeroes.  As $K$-points are Zariski dense in projective space, we can take $f$ to be defined over $K$.
Thus, $f$ defines a morphism $X\to \bbp^1$ of degree $N$, which gives a degree $N$ extension of function fields $K(X)/K(\bbp^1)$.

The Galois group $H$ of the minimal Galois extension of $K(\bbp^1)$
containing $K(X)$ is a permutation group of degree $N$, defined up to conjugation.  Local monodromy considerations at $\infty$ and $0$ show that $H$ contains an $N$-cycle and an element of order $2$ with $|S|$ $2$-cycles (and therefore at least $2$ fixed points, as $N\ge 2|S|+2$.)
If $N$ is prime and sufficiently large, $H$ cannot therefore be contained in the group of affine transformations on the field with $N$ elements.
If, in addition, $N>23$ and $N$ is not of the form $\frac{q^m-1}{q-1}$ for any prime power $q$ and integer $m\ge 2$, then $H$ must contain $A_N$ \cite[Theorem~4.2]{Feit}.  We conclude that there exists
an $n$-tuple $(h_1,\ldots,h_n)\in H^n$ such that
$$\langle h_1,\ldots,h_n\rangle = \iota(G).$$

There is no difficulty choosing an arbitrarily large prime that is not of the form $\frac{q^m-1}{q-1}$.  For the latter expression to be prime, either $q$ must be a power of $2$, or $m$ must be odd.
The number of expressions of the form $\frac{2^{km}-1}{2^k-1}$ less than $M$ is $O(\log^2 M)$; the number of expressions $1+q+\cdots+q^{2k} < M$ where $k$ is a fixed positive integer and $q\ge 2$ is variable is less than
$M^{1/2k}-1$, so the number if $k$ is allowed to vary is $O(M^{1/2})$, while the number of primes $<M$ grows like $M/\log M$.  So we may fix a prime $N$ with $A_N\subset H\subset S_N$.

Suppose $c\in K\subset  \bbp^1(K)$ such that the kernel $G_L$ of the action of $G_K$ on the $N$-element set $f^{-1}(c)$
satisfies $G_K/G_L\cong H$.  Defining
$$V_c:=\Span_{\bbq} f^{-1}(c)\subset A(L)\otimes\bbq,$$
the action of $\Gal(L/K)$ on this space is a quotient of the $N$-dimensional permutation
representation of $H$, which decomposes as an irreducible $N-1$-dimensional representation and a trivial representation.

If the quotient $V_c$ does not contain the $N-1$-dimensional irreducible,
it must have trivial $H$-action, so $x_1-x_2\in A(L)_{\tor}$ for all $x_1,x_2\in f^{-1}(c)$.  By a theorem of  Geyer and Jarden \cite[Proposition~1.1]{GJ}, since $x_1-x_2$ is defined over
a bounded degree extension of a given finitely generated field, there are only finitely many possibilities for it, independent of the choice $c$.

For each $t\in A(\bar K)_{\tor}$ there are three possibilities for the translation map $\tau_t\colon A\to A$:
\begin{enumerate}
\item $\tau_t(X)\neq X$;
\item $\tau_t(X)=X$ but $f\circ \tau_t|_X$ and $f$ are distinct;
\item $f\circ \tau_t|_X = f$
\end{enumerate}
In case (1), $f(x_1)=f(x_2)$ can only happen for finitely many pairs $(x_1,x_2)$, so there are finitely many possibilities for $c = f(x_1)$.  In case (2), there are finitely many $c$ for which there exists $x_1$ with $x_1,x_1+t\in f^{-1}(c)$.

The set of torsion points $t$ for which case (3) occurs forms a finite subgroup $T$ of $A$, and the morphism $f$ factors through $X\to X/T$.  Since $\deg f = N$ is prime, this means $|T|=1$ or
$|T| = N$.  The case $|T|=1$ means $x_1-x_2\in T$ implies $x_1=x_2$, so this case can be disregarded.  This leaves the case $|T|=N$, which implies $f$ is the quotient map by translation by $T$.  This cannot happen, since $T$ acts freely on $A$ and therefore on $X$, while every degree $N$ morphism from an irreducible curve to $\bbp^1$ is ramified.  Thus, at the cost of excluding finitely many values of $c$, we may assume that
the action of $H$ on $V_c$ contains an irreducible factor of degree $N-1$.

If the composition of $e_{\bs}$ with the quotient map $G_K\to H$ has image $\iota(G)$, then $\rho_A\circ e_{\bs}$ contains at least one copy of the atom $\alpha$.
The probability that this occurs for a single value $c$ is at least $N!^{-n}$.
However, by  Hilbert irreducibility, for every finite sequence $K_1,\ldots,K_m$ of finite extensions  of $K$ there exists $c\in K\subset \bbp^1(K)$ such that the kernel $G_L$ of the action of $G_K$ on the $N$-element set $f^{-1}(c)$
satisfies $G_K/G_L\cong H$, and $L$ is linearly disjoint from $K_1\cdots K_m$ over $K$.  We define $K_{m+1}$ to be $L$, and iterate.
By linear disjointness, the conditions on the compositions of $e_{\bs}$ with different maps $G_K\to \Gal(K_i/K)$ are independent.  By the second Borel-Cantelli lemma, this implies that
with probability $1$, the atom $\alpha$ occurs infinitely many times in $\rho_A\circ e_{\bs}$.
\end{proof}

\section{Diophantine Geometry}

Let $K$ be any Hilbertian field.
Given an abelian variety $A/K$, a finite group $G$, and an $n$-dimensional integral representation $G\to \Aut(\Lambda)$, we define $B := \Hom(\Lambda,A)\cong A^n$.  The action of $G$ on $\Lambda$ determines an action of $G$ on $B$.  Note that $\Lambda\subset \Hom(B,A)$, and this embedding is compatible with $G$-actions, where $G$ acts on $\Hom(B,A)$ through its action on $B$.

Suppose that the quotient variety $B/G$ contains a rational curve $C$.  Let $Y$ denote the inverse image of $C$ in $B$.  If $Y$ is irreducible, we can apply Hilbert irreducibility to the morphism $f\colon Y\to Y/G = C$ to conclude that for ``most'' $c\in C(K)$, the inverse image of $c$ in $Y$ consists of a single point $y\in Y$ whose residue field $L$ is a $G$-extension of $K$.  The embedding of $Y$ in $B$ gives a point of $B(L)$, which we again denote $y$.

The action of $\Gal(L/K)$ on the Galois orbit of $y$ is compatible with the action of $G$ on $B$.
Composing the embedding of $\Lambda$ in $\Hom(B,A)$ with the evaluation map on $y$,
we get a $G$-equivariant map from $\Lambda$ to $A(L)$.  The span of the image in $A(L)\otimes\bbq$, regarded as $\Gal(L/K)$-representation is therefore a quotient representation of $W := \Lambda\otimes \bbq$
as $G$-representation.
Generically \cite[Proposition~2.1]{ILPre}, this span is in fact isomorphic to $W$ as $G$-representation.

In favorable situations, this construction can be used to give unconditional proofs that $A(\bar K(\bs))$ has infinite rank.  If $W^H\neq 0$ for all subgroups of $G$ generated by $n$ elements, then a copy of $W$ in $A(L)\otimes \bbq$ guarantees a non-zero element in $A(L\cap \bar K(\bs))\otimes \bbq$.
By linear disjointness arguments, one can construct an infinite linearly independent sequence of such elements.

In \cite{Larsen}, this idea is implemented in the following concrete form.  Suppose $E$ is an elliptic curve over $K$ and $a_i$ and $b_i$ are period-$3$ sequences in $K$ such that for all $i$
(or equivalently for $i=1,2,3$),
$$y^2 = (x-a_i)(x-a_{i+1})(x-b_i)(x-b_{i+1})$$
is $K$-isomorphic to $E$.  For all $c\in K$ distinct from all $a_i$ and $b_i$, we have
$$\prod_{i=1}^3 (c-a_i)(c-a_{i+1})(c-b_i)(c-b_{i+1}) \in (K^\times)^2.$$
Thus for all $\sigma\in \Gal(\bar K/K)$ there exists $i\in \{1,2,3\}$ such that
$$\sqrt{(c-a_i)(c-a_{i+1})(c-b_i)(c-b_{i+1})}\in \bar K(\sigma)$$
Finding such $a_i$ and $b_i$ for a given $E$ amounts to realizing the $c$-line in $E^3/G$,
where $G$ is the Klein $4$-group with each non-zero element inverting two coordinates of $E^3$.

Ideally, one would like to choose $G$ and $\Lambda$ such that $\Lambda^G = (0)$ but $\Lambda^H \neq (0)$
for all $n$-generated subgroups of $G$, where $n$ is a fixed (possibly large) integer.
For any $n$, such a pair exists; for instance, one can take $G = (\bbz/2\bbz)^{n+1}$,
and let $\Lambda$ denote the quotient of the integral regular representation by its group of $G$-invariants.  Unfortunately, for such pairs $(G,\Lambda)$, for $n\ge 2$, we do not know if $\Hom(\Lambda,A)/G$ has any rational curves.  (There is one exception: for the elliptic curve
$$y^2 = (x-1)(x-\zeta_7)(x-\zeta_7^2)(x-\zeta_7^4),$$
and $n=2$, the Hamming code gives \cite[Theorem~6]{Larsen} a rational curve on $\Hom(\Lambda,A)/G$ .)

There are example of pairs $(G,\Lambda)$ for which $\Hom(\Lambda,A)/G$ has many rational curves.  For instance, a theorem of Eduard Looijenga \cite{Looijenga} shows that if $G$ is a Weyl group and $\Lambda$ is suitably chosen, this quotient is actually a weighted projective space.  On the other hand, examples where most points in  $\Hom(\Lambda,A)/G$ lie on rational curves are rare and never too far from the reflection group case \cite{KL}.
It may also happen that $\Hom(\Lambda,A)/G$ may fail to be uniruled and still have some rational curves.  There are some interesting examples of this phenomenon, especially when $A$ is an elliptic curve; for instance most of the smaller sporadic groups can be realized in this way \cite{ILPre}.

On the other hand, there is some evidence that, especially for quotients of higher dimensional abelian varieties by finite groups, rational curves may be the exception rather than the rule.
For instance, by a theorem of Gian Petro Pirola \cite{Pi}, most Kummer surfaces of dimension $\ge 3$ have no rational curves.  It would be interesting to have a criterion in terms of a positive integer $d$, a finite group $G$, and an integral representation $\Lambda$ of $G$ for whether $\Hom(\Lambda,A)/G$ has a genus $0$ curve for all abelian varieties $A$ of dimension $d$.

Note that if $\Lambda_1$ and $\Lambda_2$ are integral representations such that the rational representations $W_i := \Lambda_i\otimes\bbq$ are isomorphic to one another, then the $\Hom(\Lambda_i,A)$ admit $G$-equivariant isogenies in both directions, so the existence of genus $0$ curves depends only on the underlying rational representation.  As a very preliminary step in this direction,
in \cite{IL-age}, we gave the following sufficient criterion, extending the result of \cite{KL}:

\begin{thm}
\label{age}
If $G$ has an element $g$ which acts on the Lie algebra of $B$ with eigenvalues $\lambda_i = e^{2\pi i x_i}$, $0\le x_i < 1$, and $\sum_i x_i \le 1$, then $B/G$ has at least one rational curve.
\end{thm}

If $G$ is an alternating group of even degree $d$ and $W$ is the $d-1$-dimensional irreducible representation, then $W^H\neq 0$ for all cyclic subgroups $H$ of $G$.  Indeed, $H$ is generated by a single element, which cannot be a $d$-cycle and must therefore have at least two orbits in its action on $\{1,2,\ldots,d\}$; it follows that $W^H$ has dimension at least $1$.
Using this construction, in \cite{IL-cyclic} we proved Conjecture~\ref{conj-ii} for $n=1$:

\begin{thm}
\label{cyclic}
Let $A$ be a non-trivial abelian variety over a field $K$ which is not locally finite, which does not have characteristic $2$, and which has topologically cyclic Galois group. Then the rank of $A$ over $K$ is infinite.
\end{thm}

\section{Combinatorics}

According to Conjecture~\ref{ample}, every pointed non-singular curve $X/K$ has infinitely many points over $\bar K(\bs)$ for all $\bs$.  This implies the following:

\begin{conj}
\label{curve}
If $K$ is a finitely generated field over $\bbq$ and $X$ is a pointed non-singular curve over $K$, then for all $\bs$, there are infinitely many
$\bar K(\bs)$ points on $X$.
\end{conj}

For some curves, this can be proved directly \cite{IL1}:

\begin{thm}
If $a_1,\ldots,a_{2g+2}$ are pairwise distinct elements of any infinite field $K' := \bar K(\bs)$, not of characteristic $2$, with finitely generated Galois group, then the (non-singular) split hyperelliptic curve
$$y^2 = (x-a_1)\cdots(x-a_{2g+2})$$
has infinitely many points over $K'$.
\end{thm}

To explain the strategy of proof, we consider the case that the characteristic of $K$ is $0$, $g=1$, and $a_i = 1-i$ for $i=1,\ldots,4$.
As $G_{K'} = \langle \bs\rangle$ is finitely generated, by Kummer theory, $(K')^\times \otimes\bbf_2$ is finite.  Each positive integer $n$ determines a class in this group,
so by van der Waerden's theorem, there exist four positive integers in arithmetic progression, $a$, $a+d$, $a+2d$, $a+3d$,
all of which are equivalent modulo squares in $K'$.  Thus
$$(a/d)(a/d+1)(a/d+2)(a/d+3) \in (F')^2,$$
so there is a rational point on the curve with $x=a/d$.
In the general case, the proof uses the Hales-Jewett Theorem \cite[Theorem~1]{HJ}.

This result implies that if $K'$ is not locally finite, then  for any abelian variety $A$ which admits a non-constant $K'$-morphism from a split hyperelliptic curve $X$, the rank of $A$ over $K'$ is infinite.
In particular, this is the case for all elliptic curves with all $2$-torsion points rational.

The other class for which we are aware of a combinatorial proof of Conjecture~\ref{curve} is projective curves $X$ of the form $ax^n+by^n + cz^n = 0$.
By \cite{BL}, if $a,b,c\in \bbq$  then $X$ has infinitely many points over $\bar\bbq(\bs)$ for all $\bs$ provided it has at least one point over $\bbq$.  In fact,
less suffices; it is enough that $X$ has points over every completion of $\bbq$.  The crucial point is to solve equations of the form $au+bv+cw = 0$,
where $u,v,w$ lie in chosen cosets of an arbitrary finite index subgroup of $\bbq^\times$.  This is done using the circle method.

\section{Arithmetic}

In this section and the next, we  focus on Conjecture~\ref{weak} in the case that $A=E$ is an elliptic curve.  Here we consider methods from arithmetic geometry; in the following section, which is mainly conjectural, we consider what might be hoped for from analytic number theory.

We start with an elliptic curve $E/\bbq$.  By modularity, we have a good supply of algebraic points on $E$, namely, the Heegner points.

 Let  $K$ be an imaginary quadratic field of $\bbq$ with discriminant $D$ and $N$ be the conductor of $E/\bbq$. For each positive integer $c$ relatively prime to $ND$, let $\mathcal{O}_c$ be the order of index $c$ in the ring of integers $\mathcal{O}_K$ of $K$. Then the elliptic curve $\bbc/ \mathcal{O}_c$ defines a point on the modular curve $X_0(N)$ and its image $P_c$ on $E$ under the modular parametrization of $X_0(N)$ to $E$ is called the Heegner point of conductor $c$ and it is defined over the ring class field $H_c$ of conductor $c$.  If $K$ satisfies the so-called Heegner hypothesis, i.e., all primes dividing the conductor $N$ of $E$ split in $K$, then there exists a non-torsion Heegner point in the collection of all Heegner points. (See \cite{Dar}.)

In \cite{I2}, it is proved that Heegner points span an infinite-dimensional subspace of the Mordell-Weil group $E(H)$ over the compositum $H$ of all ring class fields with conductor prime to $ND$. In particular, since the ring class fields $H_{rp^m}$, where $r$ and $p$ are relatively prime to $ND$, have dihedral Galois group  over $\bbq$, if an automorphism $\sigma\in G_{\bbq}$ does  not fix $K$, then by the norm-compatibility relation among the Heegner points over $H_{rp^k}$, it can be shown that the rank of $E$ is unbounded over the fixed subfields of $H_{rp^m}$ under $\sigma$ as $m$ increases. If $\sigma$ fixes all imaginary quadratic extensions, then automatically the rank of $E$ over the fixed subfield under $\sigma$ is infinite by an elementary argument.
Either way, $E$ has infinite rank over $\bar\bbq(\bs)$ for $n=1$.

 The same strategy has been applied \cite{BreI} to extend  this result  to elliptic curves over global function fields of odd characteristic parametrized by Drinfeld modular curves, and for elliptic curves over totally real fields parametrized by Shimura curves.

If $k$ is a totally real number field and $f$ is a new form on $\text{GL}_2(\bba_k)$ of weight $2$ with level condition associated with $\mathfrak{c}$, where $\bba_k$ denotes the ring of ad\`{e}les of $k$ and  $\mathfrak{c}$ is a non-zero ideal of $\mathcal{O}_k$, then by \cite{Zh} there exists an elliptic curve $E'/k$ of conductor $\mathfrak{c}$ such that the $L$-functions of $E'$ and $f$ coincide up to factors at primes dividing $\mathfrak{c}$ and there exists a Shimura curve $X/k$ and a surjective $k$-morphism from $X$ to $E'$. So we say that $E/k$ has a modular parametrization by a Shimura curve if $E$ is $k$-isogenous to $E'$ arising as above. If $K$ is an imaginary quadratic extension of $k$ such that the discriminant of $K/k$ is prime to $\mathfrak{c}$ and satisfying a splitting or non-splitting property depending on the degree of $k$ over $\bbq$, which is the Heegner hypothesis in the case of totally real fields, and if $E/k$ has a modular parametrization, then we can construct Heegner points on $E$ via the isogeny as before.

If $k$ is a global function field of odd characteristic, for any elliptic curve $E/k$, there exists a morphism from the Drinfeld modular curve $X_0(\mathfrak{c})$ to $E$, where the conductor of $E$ is $\mathfrak{c}\cdot \infty$ for an ideal $\mathfrak{c}$ of $\mathcal{O}_k$ \cite{GR}. We say  $K/k$ is an \emph{imaginary quadratic extension} if the place $\infty$ does not split in $K/k$. In this case, the Heegner hypotheses for $K$ is the condition that all primes dividing $\mathfrak{c}$ split in $K/k$. For $K$ satisfying the Heegner hypothesis, we can construct Heegner points via this Drinfeld modular curve parametrization as before.

Although these results are now encompassed by Theorem~\ref{cyclic}, Heegner point methods, possibly in conjunction with ideas from analytic number theory, remain a viable approach to Conjecture~\ref{weak} for elliptic curves over suitable global ground fields.

Tim and Vladimir Dokchitser \cite{DD} gave an approach to Conjecture~\ref{weak} for elliptic curves for general $n$.
Assuming the Birch-Swinnerton-Dyer conjecture,
they proved that the conjecture holds if $K$  has a real place or $E$ has non-integral $j$-invariant.

Their main idea was to show first that there is a quadratic extension $M/K$ where the root number $w(E/M)=-1$ and the rank of $E(M)$ is odd. If a place $v$ of $K$ which is real or nonarchimedean with $v(j(E))<0$  is fixed, then they proved that by the weak approximation theorem, $M$ can be taken as a quadratic extension such that the negative contribution of the local root number occurs only above $v$, so that the global root number $w(E/M)$, which is the product of local root numbers, is $-1$.

If $G\subset G_K$ is finitely generated,  for an odd prime $p$, they constructed a Galois extension $F/K$ containing $M$ such that $\text{Gal}(F/K)\cong \bbf_p^r\rtimes C_2$ where $C_2$ acts by $-1$, the image of $G$ in $\text{Gal}(F/K)$ has order at most $2$, and the primes of $M$ above bad reduction primes of $E/K$ split completely in $F$.  Then, by taking any index $p$ subgroup $V$ of $\bbf_p^r$,  $V$ is a normal subgroup of $\text{Gal}(F/K)$ and $\text{Gal}(F^V/K)$ is isomorphic to the dihedral group of order $2p$. As the image of of $G$ in $\text{Gal}(F/K)$ has order $\le 2$, they got a degree $p$-extension $\subseteq F^V$ over $K$ fixed by $G$ and applied the congruence modulo $2$-relation among the ranks of $E$ over $M, L$ and $K$ (\cite{DD2}) from which they concluded that by the $p$-parity conjecture,
$$\rk_L E > \rk_K E.$$

\section{Analytic Number Theory}
The method of \S3 is based on the idea that where there are rational curves, there are rational points. However, many varieties without rational curves nevertheless have many rational points.
The following conjecture seems to us likely to be true but unlikely to follow from the method of \S3:

\begin{conj}
\label{paras}
Let $E/\bbq$ be an elliptic curve.  For all $n$ there exist infinitely many linearly independent $n+1$-dimensional $\bbf_2$-subspaces  $V_i\subset\bbq^\times/(\bbq^\times)^2$
such that for all $i$ and all non-zero $v\in V_i$, the twist $E_v$
has positive rank.
\end{conj}

This conjecture implies that  Conjecture~\ref{weak} holds for $E$.  Indeed, let $\bbq(\sqrt{V_i})$ denote the extension of $\bbq$ obtained by taking square roots of coset representatives of all elements of $V_i$.
By Kummer theory, $\Gal(\bbq(\sqrt{V_i})/\bbq)$ is dual to $V_i$, and the image of
$\langle \bs\rangle$ fixes  a subspace of $V_i$ of positive dimension.  Thus there exists $v_i\in V_i$ such that $\bbq(\sqrt{v_i})\subset \bar \bbq(\bs)$.
As $E_{v_i}$ has positive rank, the non-trivial eigenspace of the action of $\Gal(\bbq(\sqrt {v_i})/\bbq)$ on $E(\bbq(\sqrt{v_i}))\otimes \bbq$ has positive rank.
By taking points of $E$  in  linearly disjoint fields $\bbq(\sqrt{v_i})$, we obtain an infinite sequence of elements in $E(\bar{\bbq}(\bs))\otimes\bbq$,
which can easily be seen to be linearly independent (see, e.g., \cite[Theorem~5]{Larsen}).

It is known \cite{ILR} that there are examples of curves $E$ such that Conjecture~\ref{paras} holds for $n=1$.  The conjecture would hold in general if we knew:

\begin{conj}
\label{big-para}
Let $E/\bbq$ be an elliptic curve.  There exists an infinite dimensional $\bbf_2$-subspace $V_0\subset\bbq^\times/(\bbq^\times)^2$
such that for all non-zero $v\in V_0$, the twist $E_v$
has positive rank.
\end{conj}

Using analytic arguments and the known positive density of positive rank quadratic twists of certain elliptic curves over $\bbq$ (see, e.g. \cite{Vatsal}),
we proved  \cite{IL-para} that there exist elliptic curves $E$ over $\bbq$ with arbitrarily large finite subspaces $V_0\subset \bbq^\times/(\bbq^\times)^2$ satisfying the positive rank condition for $E$.
On the other hand, Conjecture~\ref{big-para} would be an easy consequence of the following conjecture:

\begin{conj}
\label{twists}
Let $E_1,\ldots,E_n$ be elliptic curves over $\bbq$.  Then there exists $d\in \bbq^\times/(\bbq^\times)^2$ such that the twist of every $E_i$ by $d$ has positive rank.
\end{conj}

By estimating $\sum_{d<N} L'(1,E_d)$, Perelli-Pomyka\l{}a \cite{PP} showed that the proportion of twists of a fixed elliptic curve $E$ by $d < N$ which have analytic rank $1$ grows faster than $N^{-\epsilon}$ for all $\epsilon > 0$.  The same is true for rank $1$ by Kolyvagin's theorem.
Choosing $\epsilon < 1/n$, this gives a heuristic argument for Conjecture~\ref{twists}.  It would be interesting to try to bound
$$\sum_{d<N} L'(1,(E_1)_d)\cdots L'(1,(E_n)_d)$$
away from zero.


\begin{thebibliography}{DD2}

\bibitem[BF]{BF}Bary-Soroker, Lior; Fehm, Arno:
Open problems in the theory of ample fields. Geometric and differential Galois theories, 1--11,
S\'emin. Congr., 27, Soc. Math. France, Paris, 2013.

\bibitem[BL]{BL}Bourgain, Jean; Larsen, Michael: A Finitary Hasse Principle for Diagonal Curves,  arXiv:1404.2849.

\bibitem[BI]{BreI} Breuer, Florian; Im, Bo-Hae:
Heegner points and the rank of elliptic curves over large extensions of global fields.
\textit{Canad.\ J.\ Math.}  \textbf{60}  (2008),  no.\ 3, 481--490.

\bibitem[DD1]{DD} Dokchitser, Tim; Dokchitser, Vladimir: A note on Larsen's conjecture and ranks of elliptic curves.
\textit{Bull.\ Lond.\ Math.\ Soc.}  \textbf{41} (2009), no.\ 6, 1002--1008.

\bibitem[DD2]{DD2} Dokchitser, Tim; Dokchitser, Vladimir: On the Birch-Swinnerton-Dyer quotients modulo squares.
\textit{Ann.\ of Math. (2)}  \textbf{172}  (2010),  no.\ 1, 567--596.

\bibitem[Dar]{Dar} H. Darmon, Rational points on modular elliptic curves, CBMS Regional Conference Series in Mathematics, 101. Published for the Conference Board of the Mathematical Sciences,
Washington, DC; by the AMS, Providence, RI, 2004. MR2020572 (2004k:11103)

\bibitem[FP]{FP}Fehm, Arno; Petersen, Sebastian:
On the rank of abelian varieties over ample fields.
\textit{Int.\ J.\ Number Theory} \textbf{6} (2010), no.\ 3, 579--586.

\bibitem[Fe]{Feit}Feit, Walter:
Some consequences of the classification of finite simple groups.
The Santa Cruz Conference on Finite Groups (Univ. California, Santa Cruz, Calif., 1979), pp. 175--181,
Proc.\ Sympos.\ Pure Math., 37, Amer. Math. Soc., Providence, R.I., 1980.

\bibitem[FyJ]{FJ}Frey, Gerhard; Jarden, Moshe:
Approximation theory and the rank of abelian varieties over large algebraic fields.
\textit{Proc.\ London Math.\ Soc.\ (3)} \textbf{28} (1974), 112--128.

\bibitem[FdJ]{FriedJ}Fried, Michael D.; Jarden, Moshe: Field arithmetic. Third edition. Revised by Jarden. Ergebnisse der Mathematik und ihrer Grenzgebiete. 3. Folge. A Series of Modern Surveys in Mathematics, 11. Springer-Verlag, Berlin, 2008.

\bibitem[GR]{GR} Gekeler, Ernst-Ulrich; Reversat, Marc: Jacobians of Drinfeld modular curves. \textit{J.\ Reine Angew.\ Math.} \textbf{476}(1996), 27--93.

\bibitem[GJ]{GJ}Geyer, Wulf-Dieter; Jarden, Moshe:
The rank of abelian varieties over large algebraic fields.
\textit{Arch.\ Math.\ (Basel)} \textbf{86} (2006), no.\ 3, 211--216.

\bibitem[HJ]{HJ} Hales, A.~W.; Jewett, R. I.: Regularity and positional games.
\textit{Trans.\ Amer.\ Math.\ Soc.}
\textbf{106} (1963), 222--229.

\bibitem[I1]{I1} Im, Bo-Hae: The rank of elliptic curves with rational 2-torsion points over large fields, \textit{Proc.\ Amer.\ Math. Soc.} \textbf{134} (2006), 1623--1630.

\bibitem[I2]{I2}  Im, Bo-Hae: Heegner points and Mordell-Weil groups of elliptic curves over large fields.
\textit{Trans.\ Amer.\ Math.\ Soc.}  \textbf{359}  (2007),  no.\ 12, 6143--6154.

\bibitem[IL1]{ILPre}Im, Bo-Hae; Larsen, Michael: Realizing Rational Representations in Mordell-Weil Groups, arXiv: math/0401209.

\bibitem[IL2]{IL-cyclic}Im, Bo-Hae; Larsen, Michael:
Abelian varieties over cyclic fields.
\textit{Amer.\ J.\ Math.} \textbf{130} (2008), no.\ 5, 1195--1210.

\bibitem[IL3]{IL-para}Im, Bo-Hae; Larsen, Michael:
Parallelopipeds of positive rank twists of elliptic curves.
\textit{Indiana Univ.\ Math.\ J.} \textbf{60} (2011), no.\ 1, 311--318.

\bibitem[IL4]{IL1}Im, Bo-Hae; Larsen, Michael: Some applications of the Hales-Jewett theorem to field arithmetic. \textit{Israel J.\ Math.}  \textbf{198}  (2013),  no. 1, 35--47.

\bibitem[IL5]{IL-age}Im, Bo-Hae; Larsen, Michael:
Rational curves on quotients of abelian varieties by finite groups.
\textit{Math.\ Res.\ Lett.} \textbf{22} (2015), no.\ 4, 1145--1157.

\bibitem[ILR]{ILR}Im, Bo-Hae; Lozano-Robledo, \'Alvaro:
On products of quadratic twists and ranks of elliptic curves over large fields.
\textit{J.\ Lond.\ Math.\ Soc. (2)} \textbf{79} (2009), no.\ 1, 1--14.

\bibitem[JK]{JK} Junker, Markus; Koenigsmann, Jochen: Schlanke K\"orper.
\textit{J.\ Symbolic Logic} \textbf{75} (2010), no.\ 2, 481--500.

\bibitem[KL]{KL}Koll\'ar, J\'anos; Larsen, Michael:
Quotients of Calabi-Yau varieties. Algebra, arithmetic, and geometry: in honor of Yu.\ I.\ Manin. Vol.\ II, 179--211, Progr.\ Math., 270, Birkh\"auser Boston, Inc., Boston, MA, 2009.

\bibitem[Ln]{Lang}Lang, Serge:
Fundamentals of Diophantine geometry. Springer-Verlag, New York, 1983.

\bibitem[Lr]{Larsen}Larsen, Michael:
Rank of elliptic curves over almost separably closed fields.
\textit{Bull.\ London Math.\ Soc.} \textbf{35} (2003), no.\ 6, 817--820.

\bibitem[Lo]{Looijenga}Looijenga, Eduard:
Root systems and elliptic curves.
\textit{Invent.\ Math.} \textbf{38} (1976/77), no.\ 1, 17--32.

\bibitem[PP]{PP}Perelli, A.; Pomyka\l{}a, J.: Averages of twisted elliptic L-functions.
\textit{Acta Arith.} \textbf{80} (1997), no.\ 2, 149--163.

\bibitem[Pi]{Pi}Pirola, Gian Pietro:
Curves on generic Kummer varieties.
\textit{Duke Math.\ J.} \textbf{59} (1989), no.\ 3, 701--708.


\bibitem[Va]{Vatsal}Vatsal, Vinayak:
Rank-one twists of a certain elliptic curve.
\textit{Math.\ Ann.} \textbf{311} (1998), no.\ 4, 791--794.


\bibitem[Zh]{Zh} Zhang, Shou-Wu: Heights of Heegner points on Shimura curves. \textit{Ann.\ of Math.} \textbf{153} (2001), no. 1, 27--147.


\end{thebibliography}
\end{document}